\documentclass[12pt]{article}
\usepackage{epsfig,psfrag,amsmath,amssymb,latexsym}
\usepackage{amscd }
\usepackage{amsfonts}
\usepackage{graphicx}
\newtheorem{theorem}{Theorem}[section]

\newtheorem{corollary}{Corollary}[section]

\newtheorem{definition}{Definition}[section]

\newtheorem{lemma}{Lemma}[section]

\newenvironment{proof}[1][Proof]{\textbf{#1.} }{\ \rule{0.5em}{0.5em}}
\newcommand{\N}{\mathbb{N}}
\newcommand{\R}{\mathbb{R}}
\newcommand{\Z}{\mathbb{Z}}
\numberwithin{equation}{section}
\begin{document}
\title{FLUCTUATIONS OF THE LONGEST COMMON
 SUBSEQUENCE IN THE ASYMMETRIC CASE OF 2- AND 3-LETTER ALPHABETS} 
\author{Federico Bonetto\thanks{School of Mathematics, Georgia
 Institute of Technology, Atlanta, GA 30332.\hfill\break Email:
 bonetto@math.gatech.edu}\; and Heinrich 
Matzinger\thanks{University of Bielefeld
Postfach 10 01 31, 33501 Bielefeld, Germany. 
\hfill\break Email: matzing@mathematik.uni-bielefeld.de}}

\maketitle

\abstract {We investigate the asymptotic standard deviation of the
  Longest Common Subsequence (LCS) of two independent i.i.d. sequences
  of length $n$. The first sequence is drawn from a  three letter
alphabet $\{0,1,a\}$, whilst the second sequence is binary.
  The main result of this article is that
 in this asymmetric case, the standard
  deviation of the length of the LCS
is  of order $\sqrt{n}$. This confirms
Waterman's  conjecture \cite{Waterman-estimation} for this special case. 
Our result
seems to indicate that in many other situations 
the order of the standard deviation is also $\sqrt{n}$.}

\section{Introduction}
In computational genetics and  computational linguistics one
of the basic problem is to find an optimal alignment
between two given sequences $X:=X_1\ldots X_n$ and 
$Y:=Y_1\ldots Y_n$.  This requires a scoring system
which can rank the alignments. Typically  a substitution matrix
 gives the score for each possible pair of letters.
The total score of an alignment is the sum of terms for each
aligned pair of residues, plus a usually negative term for each gap (gap penalty). 
\\[3mm]
{\footnotesize Let us look at an example. Take 
the sequences $X$ and $Y$ to be binary sequences.  Let the substitution
matrix be equal to:
 $$
\begin{array}{c|cc}
{\;}&{0}&
{1}\\\hline
0&2&1\\
1&1&3\\
\end{array}
$$
With the above matrix  we get the following scores for pairs of letter:
$$s(0,0)=2,s(0,1)=s(1,0)=1,s(1,1)=3.$$
(Here, $s(a,b)$ designates the score when we align letter $a$ with letter
$b$.) Take $X=0101$ and $Y=1100$ with the above substitution matrix
and a zero gap penalty.  The optimal alignment is:
$$
\begin{array}{cccccc}
0&1&0&1&\_&\_\\
\_&1&\_&1&0&0
\end{array}
$$ 
The above alignment gives the score
$s(1,1)+s(1,1)=3+3=6$.  This is the alignment with maximal score.
We denote by $L_n$ the maximal alignment-score.  In our example $L_n=6$.}\\[3mm]
Let $\{X_i\}_{i\in\N}$ and $\{Y_i\}_{i\in\N}$ be two ergodic processes
independent of each other.  Let $L_n$ denote the  optimal alignment score
of the two finite sequences 
$
X:=X_1\ldots X_n
$
and
$
Y:=Y_1\ldots Y_n
$. Waterman \cite{Waterman-estimation} conjectured that in many cases $VAR[L_n]$ 
is of order $n$.\\
 
Throughout this paper the substitution matrix is equal to
the identity and there is no gap penalty.
In this case, the optimal score is equal to the length of the
Longest Common
Subsequence (LCS) of $X$ and $Y$.
(A common subsequence of $X$ and $Y$ is a sequence which is a
subsequence of $X$ as well as of $Y$.)\\
We take the two sequences 
$X$ and $Y$ to be i.i.d. sequences. The letters of $X$
are drawn from the three letter alphabet $\{0,1,a\}$
and $Y$ is a binary sequence.  The main result of this article
is that $VAR[L_n]$ is of order $n$.\\[3mm]
{\footnotesize
Let $X=a11a1000$ and $Y=00110011$.   We 
remove the $a$'s from $X$ and obtain the binary sequence $X^{01}=111000$. 
The length of the LCS of $X$ and $Y$
is equal to the length of the LCS of $X^{01}$ and $Y$, since no $a$'s appear
in $Y$.
In this example a LCS is $1100$ and corresponds to the following alignment:
$$
\begin{array}{cccccccccc}
1&\_&\_&1&1&0&0&0&\_&\_\\
\_&0&0&1&1&0&0&\_&1&1
\end{array}
$$
(The letters which appear  stacked one on top of the other
correspond to the letters of the common subsequence $1100$.)
}\\[3mm]
 The reader might wonder why the case considered in the present article is
relevant.  Three letters in one sequence and two in the other might 
seem an unrealistic example. Our motivation is the following:
in any i.i.d. sequence there are finite patterns (i.e. finite words) which tend
to have  below-average expected matching scores. The  number of times
any given finite pattern occurs in $X=X_1\ldots X_n$ is roughly a
binomial variable with variance proportional to $n$.  Hence,
 the number of times we observe  a given pattern in $Y$
behaves roughly like
the number of $a$'s in $Y$. The number of $a$'s in $Y$,
 decrease the optimal score linearly.  
For a given finite pattern with low average matching score the same 
should be true.\\ 
Bonetto and Matzinger \cite{bonettosimulation} simulated 
the situation where both sequences $X$ and $Y$
are binary  i.i.d..
When the zeros and one's are equally likely
(and for $n$ around $10000$), surprisingly
the standard deviation of $L_n$ is of order
$o(n^{1/3})$.  This is similar to the behavior of the Longest
Increasing Sequence (LIS) as studied by  Baik-Deift-Johansson
\cite{BaikDeiftJohansson99} and Aldous-Diaconis \cite{Aldous99}.
 On the other hand,
Durringer,
Lember and Matzinger \cite{periodic} investigate the case where the sequence
$Y$ is periodic with a short period. 
They also find the standard deviation of $L_n$ to be of order
$o(\sqrt{n})$.  Hence a small change in the parameters of the model
can change the asymptotic behaviour of the random variable $L_n$
completely.\\
Let us mention a little bit of the history of these problems:\\ 
Using a sub additivity argument, Chv\`atal-Sankoff \cite{Sankoff1} prove that
the limit
$$
\gamma:=\lim_{n\rightarrow\infty}E[L_n]/n
$$
exists. The exact value of $\gamma$ remains however unknown.
Chv\`atal-Sankoff \cite{Sankoff1} derive upper and lower bounds
for $\gamma$, and similar upper bounds were found by Baeza-Yates,
Gavalda, Navarro and Scheihing \cite{Baeza1999} using an entropy
argument. These bounds have been improved by Deken \cite{Deken},
and subsequently by Dancik-Paterson \cite{Paterson1,Paterson2}. In
\cite{bound}, Hauser, Martinez and Matzinger  developed a Monte Carlo 
and large deviation-based
method which allows to further improve the upper bounds on
$\gamma$. Their approach can be seen as a generalization of the
method of
Dancik-Paterson.\\
For sequence with many letters, Kiwi, Loebl and Matousek,
\cite{kiwi} 
have the following interesting 
result:\\
when both sequences $X$ and $Y$ are drawn from the 
alphabet $\{1,2,\ldots,k\}$ and the letters are equiprobable, then
$\gamma\rightarrow2/\sqrt{k}$
as $k\rightarrow\infty$.\\
Waterman-Arratia \cite{watermanphase} derive a law of large deviation
for $L_n$ for fluctuations on scales larger than $\sqrt{n}$.
The order of magnitude of the deviation from the mean of $L_n$ is
unknown, and in fact it is not even known if these deviations are
larger than a power of $n$. However, using first passage percolation
methods, Alexander \cite{Alexander} proves that $E[L_n]/n$ converges
at a rate of order at least $\sqrt{\log{n}/n}$.\\

Waterman \cite{Waterman-estimation}
studies the statistical significance of the results
produced by sequence alignment methods. An important problem that
was open for decades concerns the longest increasing subsequence (LIS)
of random permutations and appears to be related to the LCS-problem.
However, it is an open question to know if solutions of the LIS-problem
can be used to study the LCS problem, see Baik-Deift-Johansson
\cite{BaikDeiftJohansson99} and Aldous-Diaconis \cite{Aldous99}.\\

Another problem related to the LCS-problem is that of comparing sequences
$X$ and $Y$ by looking for longest common words that appear both in $X$
and $Y$, and generalizations of this problem where the word does not need
to appear in exactly the same form in the two sequences. The distributions
that appear in this context have been studied by
Arratia-Gordon-Goldstein-Waterman \cite{Arratia3} and Neuhauser
\cite{Neuhauser}. A crucial role is played by the Chen-Stein Method for
the Poisson-Approximation. Arratia-Gordon-Waterman \cite{Arratia1,
Arratia2} shed some light on the relation between the Erd\"os-R\'{e}nyi
law for random coin tossing and the above mentioned problem. In
\cite{Watemanextremvalue} the same authors also developed an extreme value
theory for this problem.

For a general discussion of the relevance of string comparison for
biology and of other similar problem in computational bilogy the
reader can refer to the standard texts \cite{Krogh}\cite{Pevzner}\cite{Backofen}.

\section{Main result}

Throughout this paper $\{X_i\}_{i\in\N}$ and $\{Y_i\}_{i\in\N}$ 
are
two i.i.d. sequences which are independent of each other
and which satisfy all of the following three conditions:
\begin{enumerate}
\item{} The variables $X_i,i\in \N$, have state space $\{0,1,a\}$. 
\item{}There exists $p$,  $0<p<1$   such that
\begin{equation}
P(X_1=a)=p,\qquad P(X_1=0)=P(X_1=1)=\frac{1-p}{2}.
\end{equation}
\item{} The variables $Y_i,i\in \N$, are Bernoulli variables with
  parameter $1/2$.
\end{enumerate}
When all the three conditions above are satisfied, we say we are in 
{\it case I }.
The main result of this paper is:
\begin{theorem}
\label{mainresult}
When we are in case I, there exists $k>0$ not depending on $n$, such
that for all $n\in \N$, we have
\begin{equation}
VAR[L_n]\geq k\cdot n.
\end{equation}
\end{theorem}
There is also an upper bound for the variance
$$VAR[L_n]\leq K\cdot n$$
where $K>0$ is a constant not depending on $n$.
This upper bound follows directly from the large deviation result
for LCS of  Waterman-Arratia \cite{watermanphase}.  Let us give
 this result:
\begin{lemma}
\label{largedeviation}Assume that we are in case I, then:\\
there exists a constant $c>0$ (not depending on $n$ and $\Delta$) such
that for all $n$ large enough and all $\Delta>0$, we have that:
\begin{equation}
P\left( \left|L_n-E\left[L_n\right]\right|\geq n\Delta
   \right)\leq e^{-cn\Delta^2}
\end{equation}
\end{lemma}
Theorem \ref{mainresult} and lemma \ref{largedeviation} together
imply that the typical size of $L_n-E[L_n]$ is $o(\sqrt{n})$.
More precisely,
let $D_n:=(L_n-E[L_n])/\sqrt{n}$ denote the rescaled fluctuation of
$L_n$. Then:
\begin{theorem}\label{th}
The sequence $\{D_n\}$ is tight. Moreover, the limit of any weakly
convergent subsequence of $\{D_n\}$ is not a Dirac
measure.\end{theorem}
Theorem \ref{th} is a rather direct consequence
of theorem \ref{mainresult} and lemma \ref{largedeviation}.  
We refer the reader to \cite{periodic} for the proof .\\

\section{Proof of main theorem}

Let $N^a$ designate the numbers of $a$'s in the sequence
$X=X_1X_2\ldots X_n$.  Let $X^{01}$ designate the subsequence of $X$
consisting of all the 0's and 1's contained in $X$.  In other words,
$X^{01}$ is obtained by removing the $a$'s from the finite sequence
$X$.  Thus, $X^{01}$ is a finite sequence of i.i.d. Bernoulli
variables with parameter $1/2$ with random length.  The length of the
random binary string $X^{01}$ is equal to $(n-N^a)$. \\[3mm]
{\footnotesize
Let us illustrate this with a practical example.  For $n=6$,
assume that
 $X=011a0a$ and $Y=101011$. In this case $N^a=2$ and
$X^{01}=0110$. Obviously the $a$'s from sequence  $X$
can not be matched since $Y$ does not contain any $a$'s.
Hence, The length $L_6$ of the LCS of $X$ and $Y$ is equal
to the length of the LCS between $X^{01}$ and $Y$.
The length of the LCS is $L_6=3$. There are actual three
longest common subsequences:
 $011$, $010$ and $110$.}\\[3mm]
The main idea why $L_n$ fluctuates on the scale $\sqrt{n}$
is the following: The binomial variable $N^a$ has variance of
order $o(n)$. The variable $L_n$ tends to decrease linearly with an increase
of $N^a$ (since the $a$'s are not matched and thus constitute losses).
Hence $L_n$ should also fluctuate on the scale $\sqrt{n}$.\\ 
To prove this rigorously,
we simulate 
the variable $L_n$
in a special way.
  We first simulate a variable with same distribution as 
$N^a$. (We can call it $N^a$.)
Then we generate $X^{01}$ by using a
drop-scheme of random bits. Instead of flipping a coin independently
$n-N^a$ times in a row we generate a sequence $Z^1,Z^2,\ldots$ of
binary strings where $Z^k$ has length $k$.  $Z^{k+1}$ is obtained by
adding to $Z^k$ a random bit at a random location.\\[3mm]
{\footnotesize
 For example,
assume that we have the binary string $Z^6=00010$.  There are four
possible positions where the next bit could come:
$$
\begin{array}{c|c|c|c}
{\rm position\; 1 }&{\rm position\; 2 }&
{\rm position\; 3 }&{\rm position\; 4 }\\\hline
0x0010&
00x010&
000x10&
0001x0
\end{array}
$$
where $x$ designates the possible position of the next bit.  We assign
the same probability to each of the four above possibilities and draw
one of them at random.  We flip a fair coin, and fill the previously
chosen position with the number obtained from the fair coin. If
 the position chosen is the second one and the fair coin
gives us a $1$, then we obtain $Z^7=001010$.}\\[3mm]
  We apply this scheme
recursively on $k$ and obtain a sequence of random binary strings
$Z^1,Z^2,\ldots,Z^n$.  Let $Z^k_i$ designate that $i$-th bit of the
$k$-th string.  With that notation:
$$
Z^k=Z^k_1Z^k_2\ldots Z^k_k.
$$
 Hence, $\{Z^k_i\}_{i\leq k\leq n}$ is a triangular array
of Bernoulli variables.  Let us next define the $Z^k$'s in a formal
way: let $V_k,k\in\N$ be a sequence of i.i.d. Bernoulli variables with
parameter $1/2$.  Let $T_k,k\in\N$ be a sequence of independent
integer variables, so that $\{V_k\}_{k\in\N}$ is independent of
$\{T_k\}_{k\in\N}$.  Furthermore, for $k\in \N$, let the distribution
of $T_{k+1}$ be the
uniform distribution on the set $\{2,\ldots,k\}$, (i.e. for all
$s\in \{1,\ldots,k\}$, we have that $P(T_k=s)=1/(k-1)$.)  We define $Z^k$ 
recursively in $k$:
\begin{itemize}
\label{item}
\item{}Let $Z^2:=V_1V_2$.
\item{}Given the binary string $Z^k=Z^k_1Z^k_2\ldots Z^k_k$, we
  define
 $Z^{k+1}$:
\begin{itemize}
\item{}For all $j<T_{k+1}$, let
$$
Z^{k+1}_j:=Z^{k}_j.
$$
\item{} For $j=T_{k+1}$, let
$$
Z^{k+1}_j=V_{k+1}.
$$
\item{} For $j$, such that $T_{k+1}<j\leq k+1$,
let
$$
Z^{k+1}_j:=Z^{k}_{j-1}.
$$
\end{itemize}
\end{itemize}
(Thus $V_k$ designates the $k$-th bit added and $T_k$ designates
the position where it gets added.)

 To prove the main result of this paper, we generate a variable
having same distribution as $L_n$ using the bit-drop-scheme.
Instead of
generating the sequence $X$, we generate the triangular array
$\{Z^k_i\}_{i\leq k\leq n}$ and, independently, a random number $N^a$
with binomial distribution with parameters $p$ and $n$.  Then,
we look for the longest common subsequence of $Y$ and $Z^k$
with $k=n-N^a$.\\
More precisely, let $L_n^a(k)$ designate the length of the Longest
Common Subsequence of $Z^k$ and $Y=Y_1Y_2\ldots Y_n$.  Then:
\begin{lemma}
\label{representation}
Assume that case I holds and $Z^k$ is generated independently
of $Y$ and $N^a$, according to the mechanism described above.
 Then, $L_n$ has same distribution as
$L_n^a(n-N^a)$.
\end{lemma}
\begin{proof} For every $l,k\geq 0$ we have that
  $P(L_n=l|N^a=k)=P(L_n^a(n-k)=l)$. This gives the thesis. 
\end{proof}

\medskip

We can now explain the main idea behind the proof of Theorem
\ref{mainresult}: assume $f$ is a map with bounded slope so that
$f'(x)\geq c >0$ for all $x\in\R$.
Let $B$ be any random variable. Lemma \ref{boundforvar}
tells us, that in this case, the variance
of $f(B)$ is bounded below by $c^2\cdot VAR[B]$.
 On the other hand, the  map
$k\mapsto L_n^a(\cdot )$ is very likely to increase 
 above a linear rate larger than a
constant $k_1>0$. Hence $VAR[L_n]=VAR[L_n^a(n=N^a)]$ should be
 larger then  $k_1^2VAR[N^a]$.  The most difficult part in the
proof is showing that with high probability
 the slope of $k\mapsto L_n^a(k)$ is ``everywhere'' bounded 
below by a positive constant. 
 This problem is solved in the next section.
Let us  look at the details of the proof of Theorem \ref{mainresult}:
\medskip

\begin{lemma}
\label{boundforvar}
Let $c>0$.  Assume that $f:\R\rightarrow\R$ is a map which is
everywhere differentiable and such that for all $x\in\R$:
\begin{equation}\label{cc}
\frac{df}{dx}\geq c.
\end{equation}
Let $B$ be a random variable such that $E[|f(B)|]<+\infty$  Then:
\begin{equation}
VAR[f(B)]\geq c^2\cdot VAR[B].
\end{equation}
\end{lemma}
\begin{proof} We have that $E[B]$ and $E[f(B)]$ are finite. 
  Observe that $\lim_{x\to \pm\infty} f(x)=\pm\infty$ and $f(x)$ is
  strictly increasing so that there exists $x_0\in\R$ such that

\begin{equation}
\label{b_0}
f(x_0)=E[f(B)].
\end{equation}

By the mean value theorem, we know that there exists a map
$\delta:\R\rightarrow\R$ such that for all $x\in \R$ we have

\begin{equation}
\label{tailor}f(x)=f(x_0)+f'(\delta(x))\left(x-x_0   \right).
\end{equation}
By definition of variance and 
eqs.(\ref{b_0})(\ref{tailor}) we have: 
\begin{equation}
VAR[f(B)]=E[(f(B)-f(x_0))^2]=E[f'(\delta(B))^2\left(B-x_0\right)^2]
\end{equation}
Using eq.(\ref{cc}) we get:
\begin{equation}
\label{boundVARII}
VAR[f(B)]\geq c^2 E[(B-x_0)^2].
\end{equation}
Observe that
\begin{equation}
E[(B-x_0)^2]\geq \min_yE[(B-y)^2]=VAR[B]
\end{equation}
where we used a well known minimizing property of the variance.
This immediately gives
\begin{equation}
VAR[f(B)]\geq c^2VAR[B]
\end{equation}
which finishes this proof.
\end{proof}\\[3mm]
Typically, the (random) map $k\mapsto L^a(k)$ does not strictly
increase for every $k\in [0,n]$. But it is likely that
every order $o(\ln n)$ points, it increases by a linear quantity.
Next we define an event which guarantees that the map
$k\mapsto L^a(k)$ increases linearly on the scale $o(\ln n)$:
\begin{definition}
Let $E^n_{\rm slope}$ designate the event that $\forall i,j$, such
that $0<i<j\leq n$ and $i+ k_2\ln n \leq j$, we have:
\begin{equation}
L^a(j)-L^a(i)\geq k_1|i-j|.
\end{equation}
Here $k_1, k_2>0$ designate constants which do not depend on $n$
and which will be fixed in the proofs in sects. 4,5.
\end{definition}
The above definition gives the discrete equivalent of condition
(\ref{cc}) in the case of a discrete function. Before proceeding we
need a discrete version of Lemma \ref{boundforvar}.
\medskip
\begin{lemma}
\label{discretebound}
Let $c,m>0$ be two constants. Let $f:\Z\rightarrow \Z$ be a
non decreasing map such that:
\begin{itemize}
\item{} for all $i<j$:
\begin{equation}
\label{discretecondition1}
f(j)-f(i)\leq (j-i)
\end{equation}
\item{} for all $i,j$ such that $i+m \leq j$:
\begin{equation}
\label{discretecondition2}
f(j)-f(i)\geq c\cdot (j-i).
\end{equation}
\end{itemize}
Let $B$ be an integer random variable such that $E[|f(B)|]\leq+\infty$.
Then: 
\begin{equation}
VAR[f(B)]\geq 
c^2\left(1-\frac{2m}{c\sqrt{VAR[B]}}
\right)VAR[B].
\end{equation}
\end{lemma}
\begin{proof}  Because of conditions (\ref{discretecondition1}) and
  (\ref{discretecondition2}), we can find a continuously
  differentiable map $g:\R\rightarrow\R$ satisfying the following
  conditions:
\begin{itemize}
\item{} $g$  agrees with $f$ on every integer which is a multiple of $m$.
\item{}
 $\forall x\in\R$, we have
that 
\begin{equation}\label{small}
c\leq g'(x)\leq 1.
\end{equation}
\end{itemize}
Thus, we can apply lemma \ref{boundforvar} to $g(B)$ and  find:
\begin{equation}
\label{varofgB}
VAR[g(B)]\geq c^2\cdot VAR[B].
\end{equation}
The random variable $g(B)$ approximates $f(B)$:  
\begin{equation}
|f(B)-g(B)|\leq (1-c)\cdot m
\end{equation}
Hence,
\begin{equation}
\label{vardiff}
VAR[f(B)-g(B)]\leq m^2
\end{equation}
Since, $f(B)=g(B)+(f(B)-g(B))$, we can apply the triangular inequality
and find:
\begin{equation}
\sqrt{VAR[f(B)]}\geq \sqrt{VAR[g(B)]}-\sqrt{VAR[f(B)-g(B)]}
\end{equation}
Hence:
\begin{align*}
  VAR[f(B)]\geq VAR[g(B)]&-2\sqrt{VAR[g(B)]}\cdot\sqrt{VAR[f(B)-g(B)]}=\\
  =&VAR[g(B)]\left( 1-\frac{2\sqrt{VAR[f(B)-g(B)]}}{\sqrt{VAR[g(B)]}}
  \right) .
\end{align*} 
Applying the inequalities (\ref{varofgB}) and (\ref{vardiff}) to the last
inequality above, yields
\begin{equation}
VAR[f(B)]
\geq c^2VAR[B]\left(1-\frac{2\,m}{c\sqrt{VAR[B]}}\right)
\end{equation} 
which finishes this proof.
\end{proof}
\medskip

Let $\sigma_Z$ designate the $\sigma$-algebra of the triangular array
$Z^k_i$ and $\sigma_{YZ}$ the $\sigma$-algebra of the triangular array
$Z^k_i$ and of the $Y_i$.  Thus:
$$\label{sigmalg}
\sigma_Z:=\sigma(Z^k_i|i\leq k\leq n)\qquad 
\sigma_{YZ}:=\sigma(Z^k_i,Y_j|i\leq k\leq n,j\leq n).
$$
We are now ready for the proof of the main theorem \ref{mainresult}
of this article.

\paragraph{Proof of theorem \ref{mainresult}}

By Lemma \ref{representation} it is enough to prove that there exits
$k>0$ not depending on $n$, such that:
\begin{equation}
VAR[L^a(n-N^a)]\geq kn.
\end{equation}
Note that for any random variable $D$ and any $\sigma$-field
$\sigma$, we have
\begin{equation}
VAR[D]=VAR[\;E[D|\sigma]\;]+E[\;VAR[D|\sigma]\;].
\end{equation}
Thus, since the variance is never negative, we find that
\begin{equation}
VAR[D]\geq E[\;VAR[D|\sigma]\;].
\end{equation}
Taking $L^a(n-N^a)$ for $D$ and $\sigma_{YZ}$ for $\sigma$, we find:
\begin{equation}
\label{convarbound}
VAR[L^a(n-N^a)]\geq E[\; VAR[L^a(n-N^a)|\sigma_{YZ}]  \;]
\end{equation} 
Note that the map $L^a(\cdot)$ is $\sigma_{YZ}$-measurable.  Thus,
conditional on $\sigma_{YZ}$, $L^a(\cdot)$ becomes a non-random
increasing map.  The event $E^n_{\rm slope}$ is
$\sigma_{YZ}$-measurable.  When $E^n_{\rm slope}$ holds, then the
hypotheses of Lemma \ref{discretebound} holds for $f=L^a(\cdot)$ with
$c=k_1$ and $m=k_2\ln n$. This implies that
\begin{equation}
\label{a}
VAR[L^a(n-N^a)|\sigma_{YZ}]\geq  
\left(k_1\right)^2
\left(1-\frac{2k_2\ln n}{k_1
\sqrt{VAR[N^a|\sigma_{YZ}]}}\right)VAR[N^a|\sigma_{YZ}]
\end{equation}
Since $N^a$ is a binomial variable with parameter $p$ and $n$ and is
independent from $\sigma_{YZ}$, we have that
\begin{equation}
VAR[N^a]=VAR[N^a|\sigma_{YZ}]=np(1-p).
\end{equation}
Using the last equality with inequality (\ref{a}), we obtain:
\begin{equation}
\label{b}
VAR[L^a(n-N^a)|\sigma_{YZ}]\geq  
np(1-p)\left(k_1\right)^2
\left(1-\frac{2k_2 \ln n}{k_1\sqrt{p(1-p)n}}\right)
\end{equation}
Since, $VAR[L^a(n-N^a)|\sigma_{YZ}]$ is never negative and since
inequality \ref{b} holds, whenever $E^n_{\rm slope}$ holds, we find
\begin{eqnarray}\label{nocheineandere}
VAR[L_n]&\geq&E[\; VAR[L^a(n-N^a)|\sigma_{YZ}]  \;]\geq\cr
&\geq&n\cdot P(E^n_{\rm slope})\cdot
\left[p(1-p)\left(k_1\right)^2
\left(1-\frac{2k_2\ln n}{k_1\sqrt{p(1-p)n}}\right)
\right].
\end{eqnarray}

The expression on the right side of inequality (\ref{nocheineandere})
divided by $n$ converges to
$$
P(E^n_{\rm slope})p(1-p)\left(k_1\right)^2.
$$
We will show in Lemma \ref{slope} below that $P(E^n_{\rm
  slope})\rightarrow 1$ as $n\rightarrow \infty$.  Hence, for all $n$
big enough, $VAR[L_n]$ is larger than $np(1-p)\left(k_1\right)^2/2>0$.
This finishes the proof of theorem \ref{mainresult}.

\section{Slope of $L^a(\cdot)$}

This section is dedicated to the proof of the following lemma: 
\begin{lemma}
\label{slope}
We have that:
\begin{equation}P(E^n_{\rm slope})\rightarrow 1\end{equation}
as $n\rightarrow\infty$.
\end{lemma}
We first need a few
definitions.  A common subsequence of length $m$ of the two sequences
$Z^k$ and $Y$, can be viewed as a pair of strictly increasing
functions:
$$\left(\pi,\eta\right)$$
such that $\pi:[1,m]\rightarrow [1,k]$,
$\eta:[1,m]\rightarrow [1,n]$ and
\begin{equation}
\label{condition}
\forall i\in [1,m],\;
Z^k_{\pi(i)}=Y_{\eta(i)}.
\end{equation}

\paragraph{Definitions:}
\begin{enumerate}
\item{}Let $\pi:[1,m]\rightarrow [1,k]$
and $\eta:[1,m]\rightarrow [1,n]$ be two increasing functions. 
The pair of  $(\pi,\eta)$ is called 
{\it a pair of matching subsequences of $Z^k$ and $Y$}
iff it satisfies  condition (\ref{condition}).
\item{} Let $M^k_1$ designate the set of all pairs of matching
  subsequences of $Z^k$ and $Y$.
\item{} Let $M^k_2$ designate the set of all pairs of matching
  subsequences of $Z^k$ and $Y$ of maximal length, (i.e. of maximal length
  in the set $M^k_1$.)

\item{} Let $\leq$ indicate the natural partial order relation between
  increasing functions $\pi:[1,m]\to\N$, i.e. $\pi_1\leq\pi_2$
  iff, for every $i\in [1,m]$, $\pi_1(i)\leq\pi_2(i)$. With a slight
  abuse of notation we will indicate with $\leq$ also the partial
  order induced on the pairs of increasing function $(\pi,\eta)$,
  i.e. $(\pi_1,\eta_1) \leq (\pi_2,\eta_2)$ iff $\pi_1\leq\pi_2$ and
  $\eta_1\leq \eta_2$.
 
\item{} Let $M^k\subset M^k_2$ designate the set of all $(\pi,\eta)\in
  M^k_2$ which are minimal according to the relation $\leq$, (i.e.
  minimal in the set $M^k_2$).
\item{} Let $(\pi,\eta)$ be a pair of matching subsequences of
  length $m$ and let $i\in[0,m-1]$.  We call the quadruple
\begin{equation}
\label{match}
\left(\pi(i),\pi(i+1),\eta(i),\eta(i+1)    \right),
\end{equation}
a {\it match of $(\pi,\eta)$}. If $\eta(i)+2\leq\eta(i+1)$, we call
the match a {\it non-empty} match.  If there exists $j$, such that
$\eta(i)<j<\eta(i+1)$ and $Y_j=1$, resp. $Y_j=0$, we say that the
match {\it contains} a $1$, resp. a $0$.  We also say that the match
{\it contains} the point $j$ and call the bit $Y_j$
a {\it free bit} of the match $(\pi(i),\pi(i+1),\eta(i),\eta(i+1))$.
 Sometimes we identify the match
$(\pi(i),\pi(i+1),\eta(i),\eta(i+1))$ with the couple of binary words:
$$\left(\;Z^k_{\pi(i)}Z^k_{\pi(i)+1}\ldots Z^k_{\pi(i+1)}\;,\;
Y_{\eta(i)}Y_{\eta(i)+1}\ldots Y_{\eta(i+1)}\;\right).$$
\item{} Let $0<s<t\leq n$.  We call the integer interval
$[s,t]=\{s,s+1,\ldots,t\}$ a {\it block of $Y$}, if for all $r\in
[s,t]$ we have $Y_r=Y_s$ but $Y_{s-1}\neq Y_{s}$ and $Y_{t}\neq
Y_{t+1}$.  The cardinality $|\;[s,t]\;|=s-t+1$ is called {\it
length} of the block $[s,t]$. 
\end{enumerate}

{\footnotesize Let us give an illustrative example.  Take
$Z^6=101011$, $n=9$ and $Y=111000111$. Let $(\pi,\eta)$ be defined as
follows:
$$\pi(1)=1,\pi(2)=3,\pi(3)=4,\pi(4)=5,\pi(5)=6$$
and
$$\eta(1)=1,\eta(2)=2,\eta(3)=4,\eta(4)=7,\eta(5)=8.$$
 Then, $(\pi,\eta)$ is a pair of matching subsequences
of $Z^6$  and $Y$.  The common subsequence associated with it
is:
$$Z^6_1Z^6_3Z^6_4Z^6_5Z^6_6=Y_1Y_2Y_4Y_7Y_8=11011.$$
We represent the pair of matching subsequences
$(\pi,\eta)$ using an alignment of $Z^6$  and $Y$:
$$\begin{array}{cccccccccc}
1&0&1&\_&0&\_&\_&1&1&\_\\
1&\_&1&1&0&0&0&1&1&1
\end{array}$$
 In this example $(\pi,\eta)$ contains the four following matches:
\begin{enumerate}
\item{}
$$\begin{array}{ccc}
1&0&1\\
1&\_&1
\end{array}$$
\item{}
$$\begin{array}{ccc}
1&\_&0\\
1&1&0
\end{array}$$
\item{}
$$ \begin{array}{cccc}
0&\_&\_&1\\
0&0&0&1
\end{array}$$
\item{}
$$\begin{array}{cc}
1&1\\
1&1
\end{array}
$$
\end{enumerate}  
The first match above is empty.  The second match  contains a one. Here,
$Y_3$ is a free bit of the second match.
The third match contains two zero's: $Y_5$ and $Y_6$
are free bits of the third match. The forth match is empty.
The common subsequence $11011$ is of maximal length
(among all the common subsequences of $Z^6$ and $Y$). So,
we have that $L^a(6)=5$. Hence, $L^a(7)$ can only be equal to
$5$ or $6$.\\
 What is the probability
that $L^a(7)$ is larger by one than $L^a(6)$?
 When we generate $Z^7$ by dropping the bit $V_7$  on $Z^6$, then
there are
five positions where it can fall: 
$$
\begin{array}{c|c|c|c|c|c}
{\rm position\; 1 }&{\rm position\; 2 }&
{\rm position\; 3 }&{\rm position\; 4 }
&{\rm position\; 5 }\\\hline
1x01011&
10x1011&
101x011&
1010x11&
10101x1
\end{array}
$$
where $x$ designates the possible positions of the bit $V^7$.
 Each of these positions has same probability.  Positions $1$ and $2$
correspond to the first match.  Position $3$ corresponds to the second match.
Position $4$ correspond to the third match and position 
5 corresponds to match number four.\\
 If $V_7=1$ and the bit drops
on the match which contains a one (that is match number two
corresponding to position three,
i.e. $T_7=3$), then $L^a(7)=L^a(6)+1$.  The reason is that the bit
$V^7$ can then get matched with the free $1$-bit in match two and increase 
the score $L^a(6)$ by one.
Similarly, if $V_7=0$ and the bit $V^7$ drops
on  match number three,  the score gets increased by one, since
then $V^7$ gets matched with the ``free'' zero contained in match
number three.  Hence, when $V^7$ drops
on  match number three, the result is:
$L^a(7)=L^a(6)+1$. In general $L^a(k+1)= L^a(k)+1$, if the
bit $V_{k+1}$ drops on a match which contains a bit of the same color as
to $V_{k+1}$. (By color, we mean $0$ or $1$.) 
}\\[3mm]
From the idea of the previous example, we can get a lower bound 
for the probability that the score $L^a(k)$ increases by one. 
 The bit $V_{k+1}$ is equally likely
to be equal to one or equal to zero. So, when it
drops on a nonempty match, the score has at least $50\%$ probability
to increase. Each nonempty match corresponds to at least one position.
The bit $V^{k+1}$ has $k-1$ equally likely positions.  It follows:
for any 
 pair $(\pi,\eta)$ 
of matching subsequences
of $Z^k$  and $Y$:
\begin{equation}
\label{fundamentalineq}
P\left(\;L^a(k+1)=L^a(k)+1\;\;|\;\;Z^k,Y\;\right)\geq
\frac12\cdot
\frac{\#\;{\rm of\;nonempty\; matches\;of}\;(\pi,\eta)}
{k}
\end{equation}  
if $(\pi,\eta)$ is of maximal length.
\vspace{8mm}

Let us explain at this stage the main ideas for the proof of lemma
\ref{slope}.
We distinguish two cases depending on the value of $k$.\\
We first deal with the case $k<0.45n$. In this case it easy to show
that with large probability all the bits in $Z^k$ are matched. Let $E^n_{1k}$
be the event:
\begin{equation}
E^n_{1k}:=\left\{L_n^a(k)=k\right\}
\end{equation}
and
\begin{equation}
E^n_{1}:=\bigcap_{k=1}^{0.45n}E^n_{1,k}.
\end{equation}
Observe that we have
\begin{equation}
E^n_{1}=\left\{L_n^a(k+1)-L_n^a(k)=1,\,\,\forall k<0.45n\right\}
\end{equation}
{\it i.e.} the slope of $L_n^a(k)$ is equal to 1 for all $k<0.45n$ if
$E^n_{1}$ holds. In the next section we prove the following lemma:

\begin{lemma}\label{E1}
We have
\begin{equation}
\lim_{n\to\infty}P(E^n_1)=1.
\end{equation}
\end{lemma}
\vspace{8mm}
{\footnotesize Assume that instead of looking
for a LCS, we want to know if one sequence is contained in another.
For example for given $l\in\N$, we may be interested in 
finding out if the sequence $Z^k$ is a subsequence
of $Y_1Y_2\ldots Y_l$.  For this let $\nu(i)$ be the smallest
$l$ such that $Z^k_i$ is a subsequence of $Y_1Y_2\ldots Y_l$.
Then, $\nu(1),\nu(2),\nu(3),\ldots$ defines a renewal process.
The interarrival times $I_i=\nu(i+1)-\nu(i)$ have
geometric distribution and expectation $E[I_i]=2$.
Thus, $E[\nu(i)]=2i$ and $VAR[n]=o(n)$.  From this it follows
that if we want $Z^k$ to be with high probability a subsequence
of $Y_1Y_2\ldots Y_l$, we need to take $l$ somewhat above $2k$.
Let us give a numerical example.  Take $Z^3=001$ and
$Y=10101000111$. Then, $\nu(1)$ denotes the indices of the first
$Y_i$ equal to zero.  In this case, $\nu(1)=2$. Similarly,
$\nu(2)$ is the smallest $i\geq \nu(1)$ such that
$Y_i=Z^3_2=0$.  Here: $\nu(2)=4$.  Finally, $\nu(3)$ is the smallest
$i\geq \nu(2)$, such that $Y_3=1$, hence $\nu(3)=5$.}
\vspace{1cm}
    
\noindent Let us  next give the main ideas, 
why with high probability, the slope of
$k\mapsto L^a(k)$ is increasing linearly on the domain $[0.45n,n]$.
 We use the bit-drop scheme
to prove this: we show that
typically the random map $k\mapsto L^a(k)$ has a positive
drift $\gamma >0$. 
We define: 
\begin{equation}
E^n_{2k}:=\left\{\forall (\pi,\eta)\in M^k, \#
\hbox{ of nonempty matches of}\;(\pi,\eta)\;
\hbox{is larger than }\gamma n \right\}.
\end{equation}
When $E^n_{2k}$ holds, every  pair $(\pi,\eta)\in M^k$ 
has at least $\gamma n$ non-empty matches. 
 The proportion of non-empty matches to $k$ hence is
larger or equal to $\gamma$. Using inequality 
\ref{fundamentalineq}, it follows that  
\begin{equation}\label{abovegamma}
P\left(\;L^a(k+1)=L^a(k)+1\;\;|\;\;Z^k,Y\;\right)\geq
 0.5\cdot \gamma
\end{equation} 
when 
$E^n_{2k}$ holds. Let 
 $E^n_{2}$
be the event:
\begin{equation}
E^n_{2}:=\bigcap_{k=0.45n}^nE^n_{2k}.
\end{equation}
 Inequality \ref{abovegamma} implies, that when  $E^n_{2}$ holds, 
the map $k\mapsto L^a(k)$ has positive drift $0.5\gamma>0$
for $k\in[0.45n,n]$. By large deviation it follows, that  
with high probability  $k\mapsto L^a(k)$ has positive slope
on $[0.45n,n]$ as soon as $E^n_{2}$ holds . (See  lemma \ref{X}.)\\
It remains to explain why $E^n_2$ holds with high probability.\\
Let us first summarize the general idea:\\
We proceed by contradiction.  Assume
 all the matches of $(\pi,\eta)\in M^k_2$
were empty. Then all of the following would hold:
\begin{itemize}
\item{}
$$(\eta(1),\eta(2),\eta(3),\ldots,\eta(m))=
(\eta(1),\eta(1)+1,\eta(1)+2,\ldots,\eta(1)+m)$$
where $m$ is the length of the LCS of $Z^k$ and $Y$:
$m=L^a(k)$. 
\item{}The sequence
$$Y_{\eta(1)}Y_{\eta(2)}\ldots Y_{\eta(m)}=
Y_{\eta(1)}Y_{\eta(1)+1}\ldots Y_{\eta(1)+m}$$
is a subsequence
of 
$$Z^k_{\pi(1)}Z^k_{\pi(1)+1}\ldots Z^k_{\pi(m)}.$$
\end{itemize}
Hence we would have  two independent i.i.d. sequences of Bernoulli variables
with parameter $1/2$, where one is contained in the other as subsequence.
This implies that the sequence containing the other
must be approximately twice as long. Hence
$k$ is approximately at least twice as large as $m=L^a(k)$. Thus,
 the ratio $L^a(k)/k$ is 
close to 
$50\%$ or below.  
This  is very 
unlikely,
since it is known that the $L^a(k)/k$ is typically
above
$80\%$. This is our  contradiction.\\[3mm]
 From the previous argument it follows that with high probability
 any $(\pi,\eta)\in M^k$ contains
a non-vanishing proportion $\epsilon>0$ of free bits. (Hence,
$L^a_n(k)/\eta(L^a_n(k))\geq \epsilon$.)  We need to show
that this  proportion $\epsilon$ of free bits
 generates  sufficiently many  non-empty matches: 
the free bits should not be concentrated in a too small number
of matches.\\[3mm]
{\footnotesize  Let us go back to the numerical example on page 9
to illustrate how we count the proportion of bits that are free.
In that example, the first match of
$(\pi,\eta)$
contains no free bit.  The second match contains one free bit which is a one.
The third match contains  two free bits which are zero's.  
The forth match contains
no free bit.  The sequence $Y$ contains a total of $8$ bits which
are involved in a match of $(\pi,\eta)$. (Note that the last bit 
$Y_9$ of $Y$ is not counted 
since it is not involved in a match of $(\pi,\eta)$.) 
We have a proportion of free bits to bits involved in
matches equal to:
$$3/8=(8-5)/5=\frac{\eta(L^a_n(k))-L^a_n(k)}{\eta(L^a_n(k))}=
\frac{\eta(5)-5}{\eta(5)}.$$  The 3 free bits generate two non-empty matches.}
\\[3mm]
To prove that there are more than $\gamma n$ nonempty matches two arguments
are used:
\begin{itemize}
\item{}Any pair of matching subsequence $(\pi,\eta)$ which is minimal
according to our partial order for pairs of matches satisfies:\\
every match of $(\pi,\eta)$ can contain zero's or one's but not both at 
the same time. Hence, each match of $(\pi,\eta)\in M^k$
contains  free bits from at most one block of $Y$.
\item{}With high probability, the total number of integer points
in $[0,n]$ contained in blocks of $Y$ of length $\geq D$ is very small.
(By choosing $D$ large, we make the total number of
points contained in  blocks longer than $D$, 
much smaller than the number of free
bits.)
\end{itemize}
From the two points above, 
it follows that for $(\pi,\eta)\in M^k$, the majority
of free bits  are at most $D$ per match. This ensures that
the proportion $\epsilon$ of free bits, generates
a proportion of at least order $\epsilon/D$ non-empty matches. \\[3mm]
{\footnotesize  Let us look at an example of a pair $(\pi,\eta)$
which is of maximal length but not minimal according to our order relation
on $M^k_2$.
Take $Z^7=0101101$ and $Y=00110010111$.  Define
 the pair of matching subsequences $(\pi,\eta)$ as follows:
$$\pi(1)=1,\pi(2)=2,\pi(3)=3,\pi(4)=4,\pi(5)=5,\pi(6)=7$$
and
$$\eta(1)=1,\eta(2)=7,\eta(3)=8,\eta(4)=9,\eta(5)=10,\eta(6)=11.$$
Let us represent this pair of matching subsequences by an alignment:
$$\begin{array}{cccccccccccc}
0&\_&\_&\_&\_&\_&1&0&1&1&0&1\\
0&0&1&1&0&0&1&0&1&1&\_&1
\end{array}$$
This gives the common subsequence $010111$.  The pair
$(\pi,\eta)$ is of maximal length, but it is not minimal for our order relation
on $M^k_2$:  instead of $\eta(2)=7$, take $\eta^*(2)=3$. Let otherwise 
$\eta^*$ be equal to $\eta$. Then $(\pi,\eta^*)$ is strictly
below $(\pi,\eta)$. To construct $\eta^*$ we used the fact that
a match of $(\pi,\eta)$ contained both zero's and one's.  It is always
possible to 
find  a strictly smaller pair $(\pi,\eta^*)\in m^K_2$
when a match of $(\pi,\eta)$ contains  hero's and one's at the same time.
\\  Note that $(\pi,\eta)$ contains $5$ free bits, but only
one non-empty match.  All the free bits of $(\pi,\eta)$
are concentrated in one match. The match containing all the free bits
contains several blocks. By taking a minimal pair of matching
subsequences, this kind of situation is avoided.}\\[3mm]

Let us  look at the details of the proof  of lemma \ref{slope}.
Let $L^a_l(k)$ denote the length of the  LCS of $Z^k$ and
the sequence $Y^l:=Y_1Y_2\ldots Y_{l}$. For $Y^l$ to be entirely
contained as a subsequence in $Z^k$, one needs $k$ to be approximately twice
as long as $l$.  (We have that $Y^l$ is a subsequence of $Z^k$ iff
 $L^a_l(k)=l$.) Hence, it is unlikely that that
$Y^l$ is a subsequence of $Z^k$, when $k=2l(1-\delta)$.
(Here  $\delta>0$ is a constant not depending on $l$.)
In other words, it is unlikely that:
$$L^a_l(2l(1-\delta))\geq l.$$
Similarly, it is unlikely, that $Y^l$ 
is ``close to being a subsequence of $Z^k$'', when $k=2l(1-\delta)$:
\begin{lemma} \label{delta}There exists a function $\delta:\R\to\R$ such that
  $\lim_{\epsilon\to 0}\delta(\epsilon)=0$ and 
\begin{equation}
P\left(L^a_l\left(2l\bigl(1-\delta(\epsilon)\bigr)\right)>
l(1-\epsilon)\right)\leq Ce^{-cl}
\end{equation}
for all $l>0$ and suitable constants $c>0$ and $C>0$ not depending on $l$.
(Note that the constants $c>0$ and $C>0$ may depend on $\epsilon$.)
\end{lemma}
\medskip
We can now define:
\begin{equation}
E^n_{3l}=\left\{L_{l}^a\left(2l(1-\delta(\epsilon))\right)\leq 
(1-\epsilon)l\right\}
\end{equation}
and 
\begin{equation}
E^n_{3}:=\bigcap_{k=0.2n}^{n}E^n_{3k}
\end{equation}
where $\epsilon$ is a suitable number, to be fixed in the following, and
$\delta(\epsilon)$ is given by  Lemma \ref{delta}. It follows
 that:
\begin{corollary}\label{E3}
If $\delta(\epsilon)$ in the definition of $E^n_3$ is given by lemma
\ref{delta}, we have
\begin{equation}
\lim_{n\to\infty}P(E^n_3)=1.
\end{equation}
\end{corollary}

Typically, $L_n^a(k)$ is above $80\%\cdot k$. However, to make things easier,
we prove only that it is above $65\%\cdot k$.  We define: 
\begin{equation}
E_{4k}^n:=\left\{L_n^a(k)\geq 0.65k\right\}
\end{equation}
and
\begin{equation}
E^n_{4}:=\bigcap_{k=0.45n}^nE^n_{4,k}.
\end{equation}
The next lemma is proven in the next section:
\begin{lemma}\label{E4}
We have
\begin{equation}
\lim_{n\to\infty} P(E^n_4)=1.
\end{equation}
\end{lemma}
\medskip

Let us define the event
$E^n_{6k}$: 
\begin{equation}
E^n_{6k}:=\left\{L^a_n(k)\leq (1-\epsilon)\eta(L^a_n(k)),\; \forall
  (\pi,\eta)\in M^k\right\}
\end{equation}
and
\begin{equation}
E^n_6:=\bigcap_{k=0.45n}^nE^n_{6k}.
\end{equation}
The event $E^n_{6k}$ says that any pair of matching subsequences
$(\pi,\eta)\in M^k$ has a proportion of at least $\epsilon$
free bits.  (Note that $\eta(L^a_n(k))$ is the number of the last
bit of $Y$ involved in a match of $(\pi,\eta)$. Furthermore, $L^a_n(k)$
represents the number of bits that are ``matched'' by $(\pi,\eta)$.
 Hence, $\eta(L^a_n(k))-L^a_n(k)$ is the number of ``free'' 
bits.) 
\medskip
\begin{lemma}\label{346}Take $\epsilon>0$ small enough,
so that
\begin{equation}
\label{numbers}
\frac{50\%}{1-\delta(\epsilon)}<65\%.
\end{equation}
Then, we have that, for all $k>0.45n$,
\begin{equation}\label{3k4k6k}
E^n_{3}\cap E^n_{4k}\subset E^n_{6k}.
\end{equation}
 Thus
\begin{equation}
E^n_3\cap E^n_4\subset E^n_6.
\end{equation}
\end{lemma}

\begin{proof}Let $k\in[0.45n,n]$.
 We show that if $E^n_{6k}$ does not hold and $E^n_{3}$
holds, then $E^n_{4k}$ can not hold.  This in terms implies
 \ref{3k4k6k}.\\
Let $(\pi,\eta)\in M^k$. If $E^n_{6k}$ does not hold, than the proportion
of ``free'' bits of $(\pi,\eta)$ is below $\epsilon$.  In other words:
$$\frac{L^a_l(k)}{l}\geq 1-\epsilon$$
where $l:=\eta(L^a_n(k))$. (Note that $L^a_l(k)=L^a_n(k)$, since
$(\pi,\eta)$ is of maximal length.) It follows that
\begin{equation}
\label{I}L^a_l(k)\geq l(1-\epsilon).
\end{equation}
Now, when $E^n_{3k}$ holds, then
\begin{equation}
\label{II}L^a_l(2l(1-\delta(\epsilon)))\leq l(1-\epsilon).
\end{equation}
Comparing inequality \ref{I}, with \ref{II} and noting that
the (random) map $x\mapsto L^a_l(x)$ is increasing, yields:

$$k\geq 2l(1-\delta(\epsilon))$$
and hence
$$k\geq 2\eta(L^a_n(k))(1-\delta(\epsilon))
\geq 2L^a_n(k)(1-\delta(\epsilon)).$$
From this it follows, that:
\begin{equation}
\label{contradiction}
\frac{L^a_n(k)}{k}\leq\frac{50\%}{1-\delta(\epsilon)}<65\%
\end{equation}
where the $65\%$-bound is obtained from inequality \ref{numbers}.
Inequality \ref{contradiction} contradicts $E^n_{4k}$.
\end{proof}

\medskip
To obtain $E^n_2$ we must be sure that  the 
free bits of $Y$ do not concentrate in a small amount of
of matches of $(\pi,\eta)\in M^k$.
As explained in the example on page 12,
   any match of
$(\pi,\eta)\in M^k$ can contain  $0$'s or  $1$'s, (or nothing) but not  $0$'s
and $1$'s at the same time. This is due to the minimality respect to
the ordering $<$. In fact if $(\pi(i),\pi(i+1),\eta(i),\eta(i+1))$ is a
non empty match we must have that $Y_l\not=Y_{\eta(i+1)}$ for all
$\eta(i)<l<\eta(i+1)$. Otherwise, we could match
the bit $Z_{\pi(i +1)}$ with $Y_{l}$ instead
of $Y_{\eta(i+1)}$.  This modification would yield
a pair of matching subsequences of same length but strictly smaller
according to our order relation on $M^k_2$. Thus, all the free bits
of a match  of $(\pi,\eta)\in M^k$ are
contained in only one block of $Y$.\\
 It 
is useful to see how many bits are contained in long blocks.
Let $BLOCK^{D}$ designate the set of all blocks $[i,j]\subset
[0,n]$ of $Y$ of length at least $D$. (For the definition
of blocks see the definitions at the beginning of
this section.)  Let $N^{D}$ denote the total number
of points in the sequence $Y$ which are contained in a block of length
at least D:
\begin{equation}
N^D:=\left|\;\left\{s\in [1,n]\;|\;\exists[i,j]\in BLOCK^{D}
,s\in[i,j]\right\}\;\right|.
\end{equation}

Let $E^n_5$ designate the event:
\begin{equation}
E^n_5:=\left\{N^D\leq \epsilon n/4\right\}
\end{equation}
We will show in sec. 6 that:
\begin{lemma}\label{E5}
For every $\epsilon$ there exists $D$ such
that 
\begin{equation}
\lim_{n\to\infty} P(E^n_5)=1. 
\end{equation}
\end{lemma}
\medskip

We then have the following combinatorial fact:

\begin{lemma}
\label{123}
We have that, for all $k>0.45n$:
\begin{equation}\label{12subset3}
E^n_{4}\cap E^n_{5}\cap E^n_{6k}\subset E^n_{2k}
\end{equation}
with $\gamma=\frac{0.0425\epsilon}{D-1}$.
Thus also:
\begin{equation}
E^n_{4}\cap E^n_{5}\cap E^n_{6}\subset E^n_{2}.
\end{equation}
\end{lemma}
\begin{proof}
  We prove \ref{12subset3}. The event $E^n_{6k}$ implies that
 for each $(\pi,\eta)\in
  M^k$ there are at least $\epsilon\; \eta(L^a_n(k))$ free bits.
 We have:
\begin{equation}
\label{triv1}
 \eta(L^a_n(k))\geq L^a_n(k).
\end{equation}
  When $E^n_{4}$
holds, we have that: 
\begin{equation}
\label{triv2}L^a_n(k)\geq 0.65k.
\end{equation}
Since we take $k\geq 0.45n$, inequalities \ref{triv1} and
\ref{triv2}, together imply
that  the number of free
bits of $(\pi,\eta)\in M^k$ is at least    
$$\epsilon\; 0.65\cdot 0.45n=\epsilon\; 0.2925n.$$  By $E^n_5$,
 there are at most  $0.25\epsilon n$ bits contained in 
blocks of length $\geq D$.  Thus, there
are at least $0.0425\epsilon\cdot n$ free bits contained in blocks
of length $<D$.  Recall that every match of  $(\pi,\eta)\in M^k$ contains
free bits from only one block. Hence, every
match of $(\pi,\eta)\in M^k$ can contain at most
$D-1$ 
free bits from blocks of length $< D$.  Hence, these $\epsilon\; 0.0425n$ free bits
which are not in $N^D$, must fill at least $\epsilon\; 0.0425n/(D-1)$
 matches of  $(\pi,\eta)\in M^k$. It follows that
$(\pi,\eta)\in M^k$ has at least
$0.0425\epsilon\cdot n/(D-1)$ non-empty matches.
\end{proof}
\medskip

Lemmas \ref{346} and \ref{123} jointly imply that
that $E^n_3\cap E^n_4\cap E^n_5\subset E_2^n$.  Hence:
\begin{equation}
P(E_2^{nc})\leq P(E_3^{nc})+P(E_4^{nc})+P(E_5^{nc})
\end{equation}
where $E_x^{nc}$ denotes the complement of $E_x^{n}$.
We have that
$P(E_3^{nc})$, $P(E_4^{nc})$ and $P(E_5^{nc})$
all converge to zero when $n\rightarrow \infty$.
(This follows from Lemmas \ref{E3}, \ref{E4} and \ref{E5}.)
Hence, we have that:
\begin{equation}\label{sss}
\lim_{n\to\infty}P(E_2^n)=1.
\end{equation}

Let $\sigma_k$ denote the $\sigma$-algebra:
$$
\sigma_k:=\sigma(Z_i^k,Y_j|i\leq k,\,j\leq n).
$$
It is easy to check that $E^n_{2k}$ is $\sigma_k$-measurable.  Note
that $L^a(k+1)-L^a(k)$ is always equal to one or zero.

\begin{lemma}\label{prob}
  When $E^n_{2k}$ holds, then
\begin{equation}\label{lb}
P\left .\left(L^a(k+1)-L^a(k)=1\right|\sigma_k\right)\geq 0.5\gamma .
\end{equation}
\end{lemma}
\begin{proof} This has already been explained. (See inequality
\ref{abovegamma}).
\end{proof}

\medskip
We finally observe that
\begin{equation}
P(E^{nc}_{\rm slope})\leq  
P(E^{nc}_{\rm slope}\cap (E^{n}_{2}\cap E^n_1))+P(E^{nc}_{2})+P(E^{nc}_1).
\end{equation}
Since $P(E^{nc}_1)$ and $P(E^{nc}_{2})$ both go to zero as
$n$ goes to infinity,
we only need  to prove that
\begin{equation}
P(E^{nc}_{\rm slope}\cap (E^{n}_{2}\cap E^n_1))\to 0\;\;{\rm for}\;\;
n\rightarrow\infty,
\end{equation}
to establish lemma \ref{slope}. 
\begin{lemma}\label{X}
We have that
$$P(E^{nc}_{\rm slope}\cap (E^{n}_{2}\cap E^n_1))\to 0$$
as $n\rightarrow\infty$.
\end{lemma}

\begin{proof}We can assume that $\gamma<1$.  Define
$k_1:=0.4\gamma$, so that $k_1\leq 0.4$.
Let $$\Delta(k):=L^a_n(k+1)-L^a_n(k)$$ when $E^n_{2k}$
holds, and $\Delta(k):=1$ otherwise.  From 
eq.(\ref{lb}), it follows that:
\begin{equation}
\label{martingale}
P\left .\left(\Delta(k)=1\right|\sigma_k\right)\geq 0.5\gamma .
\end{equation}
Furthermore, $\Delta(k)$ is equal to zero or one
and $\sigma_k$-measurable.
For $k\in ]0.45n,n]$, let
$$\tilde{L}^a_n(k)=L^a_n(0.45n)+
\sum_{i=0.45n}^{k-1}\Delta(i).$$ 
For $k\in [0,0.45n]$, let $\tilde{L}_n^a(k):=L^a_n(k)$.
 Note that when $E^n_{2}$
holds, then
\begin{equation}
\label{Idontkno}
L^a(k)=\tilde{L}^a(k)
\end{equation} 
for all $k\in [0,n-1]$. Introduce the event $\tilde{E}^n_{slope}$
to be the event such that $\forall i,j$, with
$0.45n<i<j\leq n$ and $i+k_2 \ln n \leq j$, we have:
\begin{equation}
\tilde{L}^a_n(j)-\tilde{L}_n^a(i)\geq k_1|i-j|.
\end{equation}
When $E^n_1$ holds, then $L^a_n(k)$ has a slope of 
one on the domain $[0,0.45]$. Hence, the slope
condition of $E^{n}_{\rm slope}$ holds on the
domain $[0,0.45n]$, since we have $k_1\leq 0.4$.
When $E^n_2$ holds, then $L^a_n(k)$ and $\tilde{L}^a_n(k)$
are equal.  It follows that when $E^n_2$ and
$\tilde{E}^n_{slope}$ both hold, then the slope condition
of $E^n_{slope}$ is verified on the domain $[0.45n,n]$.  Hence
\begin{equation}
\label{sex}
E^n_1\cap E^n_2\cap \tilde{E}^n_{slope}= E^n_1\cap E^n_2\cap E^n_{slope}.
\end{equation}
Thus
$$P(E^{nc}_{slope}\cap E^n_1\cap E^n_2)
=P(\tilde{E}^{nc}_{slope}\cap E^n_1\cap E^n_2)
\leq P(\tilde{E}^{nc}_{slope} ).$$
It only remains to prove that $P(\tilde{E}^{nc}_{slope})$
goes to zero as $n\rightarrow \infty$.  For this we can use large deviation.
 Let 
$\tilde{E}^n_{i,j}$ be the event that
$$\tilde{L}^a_n(j)-\tilde{L}_n^a(i)\geq k_1|i-j|$$
Then
$$\tilde{E}^{n}_{slope}=\bigcap_{i,j}\tilde{E}^n_{i,j}$$
where the intersection in the last equation above is taken over
all $i,j\in [0.45n,n]$ such that $i+k_2 \ln n \leq j$.
It follows that
\begin{equation}
\label{final}
P(\tilde{E}^{nc}_{slope})\leq 
\sum_{i,j}P(\tilde{E}^{nc}_{i,j})
\end{equation}
where the last sum is taken over
all $i,j\in [0.45n,n]$ such that $i+k_2 \ln n \leq j$.
Since we took $k_1=0.4\gamma$ and because of \ref{martingale},
large deviation tells us that there exists constants $c,C>0$
such that
\begin{equation}
\label{martingale2}
P(\tilde{E}^{nc}_{i,j})\leq Ce^{-c|i-j|}
\end{equation}
for all $i,j\in\N$.
(The constants $C,c$ do not depend on $i,j$.)
Take $k_2:=3/c$. With this choice, \ref{martingale2}
becomes:
\begin{equation}
\label{martingale3}
P(\tilde{E}^{nc}_{i,j})\leq Cn^{-3}
\end{equation}
when $k_2 \ln n \leq |i-j|$.  Note that
there are less than $n^2$ terms in the sum in inequality
\ref{final}. By \ref{martingale3}, each term in the sum in inequality
\ref{final}, is less or equal to $Cn^{-3}$.  Thus
inequality \ref{final} and \ref{martingale3} together imply that
$$P(\tilde{E}^{nc}_{slope})\leq \frac{C}{n}.$$ 
\end{proof}
\section{Bounds for the probabilities.}

We report in this section several proofs of the lemmas
used in sec. 4.

\begin{lemma} for every $n$ and $\nu<0.5$ we have\label{EE1}
\begin{equation}
P(L^a_n(\nu n)=\nu n)\geq 1-e^{c(0.5-\nu)^2n}
\end{equation}
\end{lemma}
\begin{proof}
  We can build a pair of matching subsequences has follows: start from
  $Z^k_1$ and match it with the first $Y_{i_1}=Z^k_1$, then match $Z^k_2$
  with the first $Y_{i_2}=Z^k_2$ such that $i_2>i_1$. We can proceed as
  before until we reach the end of the $Z^k$ or of the $Y$. More
  precisely we can define a matching $(\pi,\eta)$ such that $\pi(i)=i$
  and $\nu(i)=\inf_{l>\nu(i-1)}\{Y_l=Z^k_i\}$ (see remark after Lemma
  \ref{E1} for an explicit example). Given $Z^k$
  and $Y$ we call $T_j$ the sequence of random variables defined by
  $T_j=\nu(j)-\nu(j-1)$. Observe that the $T_j$
  is a sequence of independent random variable all with geometric
  distribution of parameter $\frac{1}{2}$. It follows that
\begin{equation}
P(L^a_n(\nu n)=\nu n)\geq P\left(\sum_{i=0}^{\nu
    n}T_i<n\right)=P\left(\sum_{i=0}^{\nu n}T_i-\frac{1}{\nu}<0\right) 
\end{equation}
but
\begin{equation}
P\left(\sum_{i=0}^{\nu l}T_i-\frac{1}{\nu}>0\right)\leq
\inf_{s>0}E\left(e^{s\left(\sum_{i=0}^{\nu n}T_i-\frac{1}{\nu}\right)}\right)
\end{equation}
Due to the independence of the $T_i$ we have
\begin{equation}
E\left(e^{s\left(\sum_{i=0}^{\nu n}T_i-\frac{1}{\nu}\right)}\right)=
E\left(e^{s\left(T_0-\frac{1}{\nu}\right)}\right)^{\nu
  n}=\left(\frac{e^s}{2-e^s}\right)^{\nu n} e^{-n s}
\end{equation}
It is easy to check that
\begin{equation}
\inf_{s>0}\left(\frac{e^s}{2-e^s}\right)^{\nu} e^{-s}\leq e^{c(0.5-\nu)^2}
\end{equation}
for a suitable constant $c$, so that we get
\begin{equation}
P(L^a_n(\nu n)=\nu n)\geq 1-e^{c(\nu-0.5)^2n}.
\end{equation}
\medskip

\end{proof}

\medskip
\noindent\begin{proof}[Proof of lemma \ref{E1}] It follows immediately from the
  above lemma.
\end{proof}

In a very similar way we can prove that
\medskip
\begin{lemma} \label{bb2} For every $k$
\begin{equation}
P(L^a_{k}(2(1-\delta)k)=k)\leq Ce^{c\delta^2 k}
\end{equation}
\end{lemma}
\begin{proof}
Observe that the only possibility for $L_n^a(k)=k$ is that the pair of
matching subsequences constructed at the beginning of the proof of
lemma \ref{EE1} has length $k$. Using the notation of that proof we
have that
\begin{equation}
P\left(L^a_{(2-\delta)k}(k)=k\right)=P\left(\sum_{i=0}^k T_i\leq
  (2-\delta)k\right) 
\end{equation}
This quantity can be evaluated as in the previous proof to obtain the
lemma.
\end{proof}

\medskip
We can now estimate the probability of $E_{3k}^n$. 
\medskip
\noindent\begin{proof}[Proof of Lemma \ref{delta}]
  Consider a subset of $S\subset [0,l]$ containing $(1-\epsilon)l$
  points.  There are $\binom{l}{l(1-\epsilon)}$ such subset. We can
  fix the sequence $Y$ on the subset $S$. We have $2^{\epsilon l}$
  $Y$'s that agree on $S$. Calling
  $\delta(\epsilon)=\epsilon+\delta'(\epsilon)$ we have, due to Lemma
  \ref{bb2}, that the probability of matching all $Y$ in $S$ is
  bounded by $e^{-\delta'(\epsilon)^2l}$. Collecting the above estimates
  we get that

\begin{eqnarray}
P\bigl(L^a_l\left(2l\left(1-\delta(\epsilon)\right)\bigr)>
l(1-\epsilon)\right)&\leq&
2^{\epsilon l}\binom{l}{l(1-\epsilon)}e^{-c\delta'(\epsilon)^2l}\leq\cr
&\leq& Ce^{[\epsilon(\ln 2+\ln\epsilon)+
(1-\epsilon)\ln (1-\epsilon)-c\delta'(\epsilon)^2]l}
\end{eqnarray}
where we have used Stirling's formula.
Thus it is enough to chose 
\begin{equation}
\delta'(\epsilon)=\sqrt{\frac{2}{c}[\epsilon(\ln
2+\ln\epsilon)+ (1-\epsilon)\ln (1-\epsilon)]}
\end{equation} 
to obtain the lemma.
\end{proof}

\medskip

\noindent\begin{proof}[Proof of lemma \ref{E4}]
  We can divide the sequences $Z^k$ and $Y$ is subsequences of length
  10 and write $L^a_k(k)<\sum_{i=1}^{k/10}L_i$ where $L_i$ is the
  longest common subsequence
  between $Y_{10(i-1)+1}\ldots Y_{10i}$ and $Z^k_{10(i-1)+1}\ldots
  Z^k_{10i}$.  From Chvatal we know that $E(L_i)=6.97844$. From a standard
  large deviation argument we get
\begin{equation}
P\left(\sum_{i=1}^{k/10}L_i<k\left(\frac{E(L_i)}{10}-\delta\right)\right)<
\left(\inf_{s<0}E\left(e^{s\left(L_0-(0.69-\delta)
\right)}\right)\right)^{\frac{k}{10}}
\end{equation}
Calling
$p(s,\delta)=E\left(e^{s\left(L_0-(0.69-\delta)\right)}\right)$ it
easy to see that $p(s,\delta)$ is smooth in $s$, $p(0,\delta)=1$ and
$\partial_s p(0,\delta)<0$ for every $\delta>0$. This implies that
\begin{equation}
\inf_{s<0} p(s,\delta)<e^{-c(\delta)}
\end{equation}
for suitable $c(\delta)>0$. This immediately give the thesis
of the Lemma.

\end{proof}
 
\medskip
Finally we prove the  lemma \ref{E5}:
\medskip
\noindent\begin{proof}[Proof of lemma \ref{E5}]
Let $\tilde{N}^D$ be the number of 
integer points in $[0,n-D]$ which are followed by at least $D$
times the same color in the sequence $Y$.  Thus,
$\tilde{N}^D$ is the number of integer points $s\in [0,n-D]$
so that
\begin{equation}
\label{Y=}
Y_s=Y_{s+1}=\ldots =Y_{s+D}.
\end{equation}
It is easy to check that
\begin{equation}
 N^D\leq D \tilde{N}^D.
\end{equation}
Let now $\tilde Y_s$, $s\in [0,n-D]$, be equal to 1 iff \ref{Y=} holds, and
0 otherwise. We find:
\begin{equation}
\label{sumtilde}
\sum_{s=1}^n\tilde{Y}_s=\tilde{N}^D.
\end{equation}
To estimate the sum \ref{sumtilde} we can decompose
it into $D$ sub sums $\Sigma_1,\Sigma_2,\ldots,\Sigma_{D}$
where 
\begin{equation}
\Sigma_i=\sum_{\genfrac{}{}{0pt}{}{s=1,\ldots,n }{ s\; {\rm mod}\;
    D=i}} \tilde Y_s 
\end{equation}
so that 
\begin{equation}\label{sumD}
\tilde N^D=\sum_{i=1}^D \Sigma_i
\end{equation}
It is easy to see that

\begin{equation}
P\left(N^D>\frac{\epsilon}{4}n\right)\leq
P\left(\tilde{N}^D>\frac{\epsilon}{4D}n\right)\leq
D\cdot P\left(\Sigma_0>\frac{\epsilon}{4D^2}n\right)
\end{equation}
where the last inequality follows from the fact that at least one of
the addends in \ref{sumD} has to be larger than
$\frac{\epsilon}{4D^2}n$.
Now, the $Y_s$ appearing in the sub sum $\Sigma_0$ are
i.i.d. Bernoulli random variable with $P(Y_s=1)=2^{-D}$. We can
apply a large deviation argument analogous to the one used in
the previous proof and obtain
\begin{equation}
P\left(\Sigma_0>(2^{-D}+\delta)\frac{n}{D}\right)\leq e^{-c(\delta)
  \frac{n}{D}}. 
\end{equation}
with $c(\delta)>0$ for $\delta>0$.
Thus it is enough to choose $D$ such that $D2^{-D}<\frac{\epsilon}{4}$
\end{proof}

\addcontentsline{toc}{section}{
References}
\markright{References}

\bibliographystyle{plain}
\bibliography{bio}

\begin{thebibliography}{10}

\bibitem{Aldous99}
David Aldous and Persi Diaconis.
\newblock Longest increasing subsequences: from patience sorting to the
  {B}aik-{D}eift-{J}ohansson theorem.
\newblock {\em Bull.\ Amer.\ Math.\ Soc.\ (N.S.)}, 36(4):413--432, 1999.

\bibitem{Alexander}
Kenneth~S. Alexander.
\newblock The rate of convergence of the mean length of the longest common
  subsequence.
\newblock {\em Ann.\ Appl.\ Probab.}, 4(4):1074--1082, 1994.

\bibitem{Arratia3}
R.\ Arratia, L.\ Goldstein, and L.~Gordon.
\newblock Two moments suffice for {P}oisson approximations: the {C}hen-{S}tein
  method.
\newblock {\em Ann.\ Probab.}, 17(1):9--25, 1989.

\bibitem{Arratia1}
R.\ Arratia, L.\ Gordon, and M.S. Waterman.
\newblock The {E}rd{\H o}s-{R}\'enyi law in distribution, for coin tossing and
  sequence matching.
\newblock {\em Ann.\ Statist.}, 18(2):539--570, 1990.

\bibitem{Arratia2}
R.\ Arratia and M.S. Waterman.
\newblock The {E}rd{\H o}s-{R}\'enyi strong law for pattern matching with a
  given proportion of mismatches.
\newblock {\em Ann. Probab.}, 17(3):1152--1169, 1989.

\bibitem{Watemanextremvalue}
Richard Arratia, Louis Gordon, and Michael Waterman.
\newblock An extreme value theory for sequence matching.
\newblock {\em Ann.\ Statist.}, 14(3):971--993, 1986.

\bibitem{watermanphase}
Richard Arratia and Michael~S. Waterman.
\newblock A phase transition for the score in matching random sequences
  allowing deletions.
\newblock {\em Ann.\ Appl.\ Probab.}, 4(1):200--225, 1994.

\bibitem{Baeza1999}
R.A.\ Baeza-Yates, R.\ Gavald{\`a}, G.\ Navarro, and R.~Scheihing.
\newblock Bounding the expected length of longest common subsequences and
  forests.
\newblock {\em Theory Comput.\ Syst.}, 32(4):435--452, 1999.

\bibitem{BaikDeiftJohansson99}
Jinho Baik, Percy Deift, and Kurt Johansson.
\newblock On the distribution of the length of the longest increasing
  subsequence of random permutations.
\newblock {\em J. Amer. Math. Soc.}, 12(4):1119--1178, 1999.

\bibitem{bonettosimulation}
Federico Bonetto and Heinrich Matzinger.
\newblock Simulations for the fluctuations of the lcs-problem.
\newblock {\em in preparation}, 2004.

\bibitem{Sankoff1}
V{\'a}cl{\'a}v Chvatal and David Sankoff.
\newblock Longest common subsequences of two random sequences.
\newblock {\em J.\ Appl.\ Probability}, 12:306--315, 1975.

\bibitem{Backofen}
Peter Clote and Rolf Backofen.
\newblock {\em Computational molecular biology}.
\newblock John Wiley, Chichester, NY, 2000.

\bibitem{Paterson1}
Vlado Dancik and Mike Paterson.
\newblock Upper bounds for the expected length of a longest common subsequence
  of two binary sequences.
\newblock {\em Random Structures Algorithms}, 6(4):449--458, 1995.

\bibitem{Deken}
Joseph~G. Deken.
\newblock Some limit results for longest common subsequences.
\newblock {\em Discrete Math.}, 26(1):17--31, 1979.

\bibitem{Krogh}
Richard Durbin, Sean Eddy, Anders Krogh, and Graeme Mitchison.
\newblock {\em Bilogical Sequence Analysis}.
\newblock Cambridge University Press, Cambridge, 1998.

\bibitem{bound}
Rafael Hauser, Servet Martinez, and Heinrich Matzinger.
\newblock Large deviation based upper bounds for the lcs-problem.
\newblock {\em Submitted}, 2004.

\bibitem{kiwi}
Marcos Kiwi, Martin Loebl, and Jiri Matousek.
\newblock Expected length of the longset common subsequence for large
  alphabets.
\newblock {\em preprint}, 2003.

\bibitem{periodic}
Jyri Lember, Heinrich Matzinger, and Clement Durringer.
\newblock Deviation from mean in sequence comparison with a periodic sequence.
\newblock {\em in preparation}, 2004.

\bibitem{Neuhauser}
Claudia Neuhauser.
\newblock A {P}oisson approximation for sequence comparisons with insertions
  and deletions.
\newblock {\em Ann.\ Statist.}, 22(3):1603--1629, 1994.

\bibitem{Paterson2}
Mike Paterson and Vlado Dancik.
\newblock Longest common subsequences.
\newblock In {\em Mathematical foundations of computer science 1994 (Kosice,
  1994)}, volume 841 of {\em Lecture Notes in Comput.\ Sci.}, pages 127--142.
  Springer, Berlin, 1994.

\bibitem{Pevzner}
Pavel Pevzner.
\newblock {\em Computational molecular biology}.
\newblock MIT Press, Cambridge, MT, 2000.

\bibitem{Waterman-estimation}
Michael~S. Waterman.
\newblock Estimating statistical significance of sequence alignments.
\newblock {\em Phil.\ Trans.\ R.\ Soc.\ Lond.\ B}, 344:383--390, 1994.

\end{thebibliography}
\end{document}